\documentclass{amsart}

\usepackage{amscd}
\usepackage{amssymb}
\usepackage{amsxtra}
\usepackage{longtable}
\usepackage{calc}

\allowdisplaybreaks[2]

\theoremstyle{plain}
\newtheorem{thm}{Theorem}
\newtheorem{cor}{Corollary}
\newtheorem{prop}{Proposition}
\newtheorem{lem}{Lemma}

\theoremstyle{definition}
\newtheorem{defn}{Definition}

\theoremstyle{remark}
\newtheorem{rmk}{Remark}
\newtheorem{rmks}{Remarks}
\newtheorem{exmp}{Example}
\newtheorem{exmps}{Examples}

\renewcommand{\epsilon}{\varepsilon}
\renewcommand{\kappa}{\varkappa}
\renewcommand{\theta}{\vartheta}

\newcommand{\IQ}{\ensuremath{\mathbb Q}}
\newcommand{\IP}{\ensuremath{\mathbb P}}

\newcommand{\IZ}{\ensuremath{\mathbb Z}}

\newcommand{\td}{\ensuremath{\mathrm{td}}}

\DeclareMathOperator{\Pic}{Pic}

\newcommand{\varcard}[1]{\ensuremath{\operatorname{card}}}
\newcommand{\ord}[1]{\ensuremath{\operatorname{ord}}}

\newcommand{\diff}{\ensuremath{\mathrm d}}

\DeclareMathOperator{\id}{id}

\newcommand{\OX}{\ensuremath{{\mathcal O_X}}}

\def\Z{\IZ}
\def\Q{\IQ}

\title{On the Chern numbers of the Generalised Kummer Varieties}
\author[M.\ A.\ Nieper]{Marc A.\ Nieper-Wi\ss kirchen}
\address{Mathematisches Institut der Univ.\ zu K\"oln \\
  Weyertal 86--90 \\ 50931 K\"oln \\ Germany}
\email{mnieper@mi.uni-koeln.de}
\thanks{I would like to thank Michael Britze and Daniel Huybrechts who
  discussed the ideas presented in this note with me. Without them it would
  have taken much more time to complete this paper.}
\date{\today}

\begin{document}

\begin{abstract}
  Let $A^{[[n]]}$ denote the $2(n - 1)$-dimensional generalised Kummer variety
  constructed from the abelian surface $A$. Further, let $X$ be an arbitrary
  smooth
  projective surface with
  $\int_X c_1(X)^2 \neq 0$, and $X^{[k]}$ the Hilbert scheme of
  zero-dimensional subschemes of $X$ of length $k$. We give a formula which
  expresses the value of any complex genus on $A^{[[n]]}$ in terms of Chern
  numbers of the varieties $X^{[k]}$.
  
  In~\cite{ellings87} and~\cite{ellings96} it is shown how to use Bott's
  residue formula to effectively calculate the Chern numbers of the Hilbert
  schemes $(\IP^2)^{[k]}$ of points on the projective plane. Since
  $\int_{\IP^2} c_1(\IP^2)^2 = 9 \neq 0$ we can use these numbers and our
  formula to calculate the Chern numbers of the generalised Kummer varieties.

  A table with all Chern numbers of the generalised
  Kummer varieties $A^{[[n]]}$ for $n \leq 8$ is included.
\end{abstract}

\maketitle

\section{Introduction}

The two main series of irreducible holomorphic symplectic complex manifolds
are Hilbert schemes of points on a K3 surface and generalised Kummer
varieties invented by Beauville~\cite{beauville83}. One can ask for their
complex cobordism class tensored with $\IQ$, which is given by the values of
all their Chern numbers.

In~\cite{lehn01} Ellingsrud, G\"ottsche and Lehn proved that the complex
cobordism class of a Hilbert scheme $X^{[n]}$ of zero-dimensional subschemes
of length $n$ on a smooth projective surface $X$ over the complex numbers
depends only on the cobordism class of the surface $X$, i.e.\ on $c_1(X)^2$
and $c_2(X)$ (here and later on, top intersections on surfaces are to be
understood as intersection numbers). They showed how this result can be used
to calculate the Chern numbers of any such Hilbert scheme $X^{[n]}$ if one
knows the Chern numbers of the varieties $(\IP^2)^{[k]}$ and $(\IP \times
\IP)^{[k]}$, which in turn can be calculated by means of Bott's residue
formula.

Therefore, the Chern numbers of the Hilbert schemes of points on a K3 surface
can be efficiently calculated though now explicit formula is known. These
numbers can, for example, be used to check the conjecture
of~\cite{dijkgraaf97} about the elliptic genus of the Hilbert schemes
$X^{[k]}$ where $X$ is a K3 surface (see~\cite{lehn01}). In this note, we want
to give a
method for computing the Chern numbers of the $(n - 1)$-dimensional
generalised Kummer variety $A^{[[n]]}$ for an abelian surface $A$ and for
general $n$ as this has not appeared in the literature so far.

The $\chi_y$-genus of $A^{[[n]]}$ has been calculated by G\"ottsche and
Soergel~\cite{goettsche93}. Expressing this genus in terms of Chern numbers by
using the Hirzebruch-Riemann-Roch formula gives us enough information to
deduce the Chern numbers of $A^{[[n]]}$ for $n \leq 4$.

Using the theory of Rozansky-Witten invariants, Sawon~\cite{sawon99} produced
a further relation that allowed him to compute all the Chern numbers for
$n \leq 5$. The Chern numbers of $A^{[[6]]}$ were calculated
by M.~Britze and the author in~\cite{britze01}. J.~Sawon informed us that he
also
had computed these numbers. However, all these methods are not sufficient
to compute the Chern numbers for $n > 6$.

We will give a closed formula describing the value of any complex genus on a
generalised Kummer variety in terms of genera of the Hilbert schemes of points
on a fixed surface $X$ with $c_1(X)^2 \neq 0$. Since these are computable by
means of Bott's residue formula
for $X$ being the projective plane $\IP^2$ we can compute the
Chern numbers of $A^{[[n]]}$ for any $n$. We have done this for $n \leq 8$
(see appendix).

\section{The generalised Kummer varieties}

Let $X$ be a smooth projective surface over the field of complex numbers. For
every nonnegative integer we denote by $X^{[n]}$ the Hilbert scheme of
zero-dimensional subschemes of $X$ of length $n$. By a result of
Fogarty (\cite{fogarty68}), this scheme is smooth and projective of dimension
$2n$. It can be viewed as a resolution $\rho: X^{[n]} \to X^{(n)}$ of
the $n$-fold symmetric product $X^{(n)} := X^n/\mathfrak{S}_n$ of $X$.
The morphism
$\rho$, sending closed points, i.e.\ subschemes of $X$, to their support
counting multiplicities, is called the Hilbert-Chow morphism.

Let us briefly recall the construction of the generalised Kummer
varieties introduced by Beauville~\cite{beauville83}.
Let $A$ be an abelian surface and $n > 0$. There is an obvious summation
morphism $A^{(n)} \to A$. We denote its composition with the Hilbert-Chow
morphism $\rho: A^{[n]} \to A^{(n)}$ by $\sigma: A^{[n]} \to A$.
\begin{defn}
 The \emph{$n^{\textrm{th}}$ generalised Kummer variety $A^{[[n]]}$}
 is the fibre of $\sigma$ over $0 \in A$.
\end{defn}

Beauville showed among other things the following property of the varieties
$A^{[[n]]}$:
\begin{prop}
  The $n^{\mathrm{th}}$ generalised Kummer variety is a smooth projective
  holomorphic symplectic variety of dimension $2(n - 1)$.
\end{prop}

\begin{proof}
  \cite{beauville83}.
\end{proof}

Since $A$ acts on itself by translation there is also an induced operation of
$A$ on the Hilbert schemes $A^{[n]}$. Let us denote the restriction of this
operation to the generalised Kummer variety by $\nu: A \times A^{[[n]]} \to
A^{[n]}$. The following diagram is a cartesian one:
\begin{equation}
  \begin{CD}
    A \times A^{[[n]]} @>{\nu}>> A^{[n]} \\
    @V{\pi_A}VV @VV{\sigma}V \\
    A @>>n> A.
  \end{CD}
\end{equation}
Here, $n: A \to A, a \mapsto na$ is the (multiplication by $n$)-morphism.
It is a Galois covering of degree $n^4$. Therefore, also $\nu$ is a
Galois covering of degree $n^4$.

Next, we want to introduce certain line bundles on the Hilbert schemes and
generalised Kummer varieties that are constructed from line bundles on the
underlying surface:

Each line bundle $L$ on a smooth projective surface $X$ gives us a line bundle
$L_n$ on $X^{[n]}$ in the following way: $L^{\boxtimes n} :=
\underbrace{L \boxtimes \dots \boxtimes L}_n$ is a $\mathfrak S_n$-invariant
line bundle on the $n^{\textrm{th}}$-product $X^n$ of $X$. Therefore, we can
define the sheaf $L^{(n)} := (\pi_* (L^{\boxtimes n}))^{\mathfrak S_n}$ of
$\mathfrak{S_n}$-invariant sections of $\pi_* (L^{\boxtimes n})$ on $X^{(n)}$
where $\pi: X^n \to X^{(n)}$ is the canonical projection. The pull-back $L_n
:= \rho^* L^{(n)}$ by the Hilbert-Chow morphism is a line bundle on $X^{[n]}$.
Note that $\Pic(X) \to \Pic(X^{[n]}), L \mapsto
L_n$ is a homomorphism of groups.

This construction has already appeared for example in~\cite{lehn01}
and~\cite{britze01}.  If $X$ is an abelian surface, we denote by $L^{[[n]]}$
the restriction of $L_n$ to the generalised Kummer variety $X^{[[n]]}
\subseteq X^{[n]}$. By using the seesaw principle (cf.~\cite{mumford70}), it
can be shown that
\begin{equation}
  \label{equ:britze}
  \nu^* L_n = L^n \boxtimes L^{[[n]]}
\end{equation}
(cf.~\cite{britze01}).

\section{Complex genera in general}

Let $\Omega := \Omega^{\mathrm{U}} \otimes \IQ$ denote the complex cobordism
ring. A complex genus $\phi$ is a ring homomorphism $\phi: \Omega \to R$ into
any $\Q$-algebra $R$. It is result of Milnor (\cite{milnor60}) that the complex
cobordism ring is a polynomial ring in the cobordism classes of the complex
projective spaces.
The $R$-valued complex genera are in one-to-one correspondence
with the formal power series $f_\phi \in R[[x]]$ over $R$ with constant
coefficient $1$. The correspondence is given as follows:
\begin{equation}
  \phi(X) = \int_X \prod_{i = 1}^n f_\phi(\gamma_i),
\end{equation}
where $X$ is any complex manifold of dimension $n$ and $\gamma_1, \dots,
\gamma_n$ are the Chern roots of $X$. Therefore, the cobordism class of a
manifold is determined by the values of its Chern numbers (take $\phi =
\id_\Omega: \Omega \to \Omega$). By Milnor's result, the converse is
also true, i.e.\ the cobordism class determines the Chern numbers.

Now, let us slightly generalise the notion of a genus.
\begin{defn}
  Let $\phi$ be a complex genus. For a complex manifold $X$ together with a
  line bundle $L$ on $X$ we define
  \begin{equation}
    \phi(X, L) := \int_X e^{c_1(L)} \prod_{i = 1}^n f_\phi(\gamma_i)
  \end{equation}
  as \emph{the genus $\phi$ of the pair $(X, L)$}.
\end{defn}

\begin{rmk}
  Obviously, $\phi(X, \OX) = \phi(X)$.
\end{rmk}

\begin{exmp}
  \label{exmp:gentodd}
  If $\td(X)$ denotes the Todd genus of $X$, and $\chi(X, L)$ the holomorphic
  Euler characteristic of the line bundle $L$ on $X$, we have by the
  Hirzebruch-Riemann-Roch theorem that
  \begin{equation}
    \td(X, L) = \chi(X, L).
  \end{equation}
\end{exmp}

The genera of pairs $(X, L)$ have the following properties, which follow
directly from the appropriate properties of Chern classes/roots.
\begin{prop}
  \label{prop:gengen}
  Let $\phi: \Omega \to R$ be any complex genus with values in $R$. We have
  \begin{enumerate}
  \item
    $\phi(X \times Y, L \boxtimes M) = \phi(X, L) \phi(Y, M)$ for two
    complex manifolds $X$ and $Y$ together with a line bundle $L$ resp.\ $M$.
  \item
    $\phi(X, \nu^* L) = \deg(\nu) \phi(Y, L)$ for any Galois covering $\nu: X
    \to Y$ and any line bundle $L$ on $Y$.
  \end{enumerate}
\end{prop}

Furthermore, any genus gives us a deformed genus in the following sense:
\begin{defn}
  Let $\phi$ be a complex genus with values in the $\Q$-algebra $R$. By
  $\phi_t$ we denote the genus with values in $R[t]$ given by
  \begin{equation}
    \phi_t(X) := \int_X \prod_{i = 1}^n \left(f_\phi(\gamma_i)
    e^{t \gamma_i}\right)
  \end{equation}
  for any complex manifold $X$.
\end{defn}

\begin{rmk}
  We have $\phi_n(X) = \phi(X, K_X^{-n})$ for any integer $n$ where $K_X$ is
  the canonical line bundle on $X$.
\end{rmk}

\section{Complex genera of Hilbert schemes of points on surfaces}

In this section we want to cite some of the results of~\cite{lehn01} and give
some corollaries which will be used later on.

Let $X$ be a smooth projective surface. Following~\cite{lehn01}, we define
\begin{equation}
  H_X := \sum_{n = 0}^\infty [X^{[n]}] z^n
\end{equation}
as an invertible element in the formal power series ring
$\Omega[[z]]$. Analogously we define
\begin{equation}
  K := \sum_{n = 1}^\infty [A^{[[n]]}] z^n
\end{equation}
in $\Omega[[z]]$ where $A$ is any abelian surface. The cobordism class does
not depend on the choice of $A$ since the generalised Kummer varieties deform
with $A$. We can reformulate our
task to determine the Chern numbers of the generalised Kummer varieties by
asking: What is the value $\phi(K) \in R[[z]]$ for any complex genus $\phi:
\Omega \to R$.

There are various tautological bundles on $X^{[n]}$ (cf.~\cite{lehn99}).
The construction is as follows: Since $X^{[n]}$ represents a functor there is
a universal family $\Xi_n \subseteq X^{[n]} \times X$. Let us denote by
$\mathcal O_n$ its structure sheaf. For any locally free sheaf $F$ on $X$ we
define the sheaf $F^{[n]} := p_*(\mathcal O_n \otimes q^* F)$ on $X^{[n]}$,
where $p: X^{[n]} \times X \to X^{[n]}$ and $q: X^{[n]} \times X \to X$ are
the canonical projections.

The following lemma is a generalization of Theorem~4.2 in~\cite{lehn01} for
line bundles.
\begin{lem}
  Let $k$ be a nonnegative integer, $m_1, \dots, m_k \in \Z$, and
  $\phi: \Omega \to R$ be a
  genus. Then there exist uniquely determined universal power series
  $A_{i, j} \in R[[z]], 1 \leq i \leq j \leq k$,
  and $B_1, \dots, B_k \in R[[z]]$, and $C, D \in R[[z]]$ depending only on
  $\phi$ and $m_1, \dots, m_k$ such that
  for every smooth projective surface $X$ and line bundles $L_1, \dots, L_k$
  on $X$ we have
  \begin{multline}
    \sum_{n = 0}^\infty \phi\left(X^{[n]}, \det(L_1^{[n]})^{m_1} \otimes \dots
      \otimes \det(L_k^{[n]})^{m_k} \right) z^n
    \\
    = \exp\left(\sum_{1 \leq i \leq j \leq k} c_1(L_i) c_1(L_j) A_{ij} +
      \sum_{i = 1}^k c_1(L_i) c_1(X) B_i + c_1(X)^2 C + c_2(X) D\right).
  \end{multline}
\end{lem}

\begin{proof}
  First note that for $k = 1$ the statement of the theorem is just Theorem~4.2
  of~\cite{lehn01} for the case of line bundles with $\Psi$ (in the notation
  of~\cite{lehn01}) being the Chern
  character of the $m_1^{\mathrm{th}}$ power of the determinant.

  Theorem~4.2 of Ellingsrud, G\"ottsche and Lehn and the proof presented
  by them can be easily generalised for more than one bundle, i.e. for $k >
  1$. Therefore, our lemma as a specialization of this generalization is
  proven.
\end{proof}

From the lemma we conclude the following:
\begin{prop}
  \label{prop:hilbgengen}
  Let $\phi: \Omega \to R$ be a genus. Then there exist uniquely determined
  universal power series $A_\phi, B_\phi, C_\phi, D_\phi \in R[[z]]$
  depending only on $\phi$ such
  that for every smooth projective surface $X$ together with a line bundle $L$
  on it we have
  \begin{multline}
    \phi(H_{X, L}) := \sum_{n = 0}^\infty \phi\left(X^{[n]}, L_n\right) z^n
    \\
    = \exp\left(c_1(L)^2 A_\phi + c_1(L) c_1(X) B_\phi + c_1(X)^2 C_\phi
      + c_2(X) D_\phi\right).
  \end{multline}
\end{prop}

\begin{proof}
  As noted in section~5 of~\cite{lehn01}, we have
  \begin{equation}
    L_n = \det(L)_n = \det(L^{[n]}) \otimes \det(\mathcal
    O_X^{[n]})^{-1}.
  \end{equation}
  Therefore by the previous lemma,
  \begin{multline}
    \phi(H_{X, L}) = \sum_{n = 0}^\infty \phi\left(X^{[n]}, \det(L^{[n]})
       \otimes \det(\mathcal O_X^{[n]})^{-1}\right)
     \\
     = \exp(
     c_1(L)^2 A_{11} + c_1(L) c_1(\mathcal O_X) A_{12}
     + c_1(\mathcal O_X)^2 A_{22}
     \\
     + c_1(L) c_1(X) B_1 + c_1(\mathcal O_X)
     c_1(X) B_2 + c_1(X)^2 C + c_2(X) D)
  \end{multline}
  for certain power series $A_{i,j}, B_i, C, D$ independent of $X$ and $L$.
  Since $c_1(\mathcal O_X) = 0$ this proves the proposition with $A_\phi =
  A_{11}, B_\phi = B_1, C_\phi = C$ and $D_\phi = D$.
\end{proof}

It is possible to express the power series $A_\phi$ in terms of genera of
Hilbert schemes of points on surfaces:
\begin{prop}
  \label{prop:aphi}
  Let $\phi: \Omega \to R$ be any genus. For every smooth projective surface
  $X$,
  \begin{equation}
    c_1(X)^2 A_\phi = \frac 1 2 \ln \frac{\phi_1(H_X)
    \phi_{-1}(H_X)}{\phi(H_X)^2}.
  \end{equation}
\end{prop}

\begin{proof}
  In~\cite{hl97} it is proven that the canonical bundle of $X^{[n]}$ is $K_n$
  where $K$ denotes the canonical bundle on $X$.

  Using this we have by proposition~\ref{prop:hilbgengen} that
  \begin{multline}
    \ln \phi_m(H_X) = \ln \sum_{n = 0}^\infty \phi_m(X^{[n]}) z^n
    = \ln \sum_{n = 0}^\infty \phi(X^{[n]}, K_n^{-m}) z^n
    \\
    = m^2 A_\phi c_1(K)^2 - m B_\phi c_1(K) c_1(X) + C_\phi c_1(X)^2
    + D_\phi c_2(X)
    \\
    = (m^2 A_\phi + m B_\phi + C_\phi) c_1(X)^2 + D_\phi c_2(X),
  \end{multline}
  for all integers $m$, which proves the proposition.
\end{proof}

\section{Complex genera of the generalised Kummer varieties}

In this section we will relate the (generalised) complex genera of Beauville's
generalised Kummer variety to the complex genera of Hilbert schemes of
points on surfaces, which we studied in the previous section.

The first step in this direction is the following:
\begin{prop}
  \label{prop:kummergengen}
  Let $\phi: \Omega \to R$ be a complex genus with values in the $\Q$-algebra
  $R$. For every abelian surface $A$ together with a line bundle $L$ on it we
  have
  \begin{equation}
    c_1(L)^2 \phi(A^{[[n]]}, L^{[[n]]}) = 2 n^2 \phi(A^{[n]}, L_n)
  \end{equation}
  for all positive integers $n$.
\end{prop}

\begin{proof}
  We will make use of~\eqref{equ:britze}. Recall that $\nu$ is a
  Galois covering of degree $n^4$. By proposition~\ref{prop:gengen} we have
  \begin{multline}
    \phi(A, L^n) \phi(A^{[[n]]}, L^{[[n]]})
    = \phi(A \times A^{[[n]]}, L^n \boxtimes L^{[[n]]})
    \\
    = \phi(A \times A^{[[n]]}, \nu^* L_n)
    = n^4 \phi(A^{[n]}, L_n),
  \end{multline}
  which proves the theorem once we have shown that $\phi(A, L^n) = \frac{n^2}
  2 c_1(L)^2$.
  This follows from the fact that the Chern classes of an abelian surface are
  trivial:
  \begin{equation}
    \phi(A, L^n) = \int_A f_\phi(\gamma_1) f_\phi(\gamma_2) e^{c_1(L^n)}
    = \int_A \frac{c_1(L^n)^2} 2 = \frac{n^2} 2 c_1(L)^2,
  \end{equation}
  where we used that $f_\phi$ is a power series with constant coefficient
  $1$.
\end{proof}

In~\cite{britze01} M.~Britze and the author expressed the (holomorphic) Euler
characteristic of the line bundle $L^{[[n]]}$ in terms of the Euler
characteristic of $L$ in order to deduce a formula for the Euler
characteristic of an arbitrary line bundle $M$ on $A^{[[n]]}$ as a polynomial
in the Beauville-Bogomolov quadratic form of $c_1(M)$. By using the analogous
expression of the Euler characteristic of the line bundle $L_n$ on $A^{[n]}$
(see~\cite{lehn01}) we get the mentioned result of~\cite{britze01} as a
corollary of the previous theorem:
\begin{cor}[\cite{britze01}]
  The holomorphic Euler characteristic of the line bundle $L^{[[n]]}$ on
  $A^{[[n]]}$ is given by
  \begin{equation}
    \chi(A^{[[n]]}, L^{[[n]]}) = n \binom{\chi(A, L) + n - 1}{n - 1}.
  \end{equation}
\end{cor}

\begin{proof}
  By lemma~5.1 of~\cite{lehn01} we have
  \begin{equation}
    \chi(A^{[n]}, L_n) = \binom{\chi(A, L) + n - 1} n.
  \end{equation}
  Using this, the corollary follows from the proposition applied to the case
  for $\phi$ being the Todd genus (remember example~\ref{exmp:gentodd}). Also
  note that $\chi(A, L) = \frac 1 2 c_1(L)^2$ by the Hirzebruch-Riemann-Roch
  formula.
\end{proof}

If we are interested in the usual genera of the generalised Kummer varieties,
i.e.\ the genera of the pairs $(A^{[[n]]}, \mathcal O_{A^{[[n]]}})$, we can't
use proposition~\ref{prop:kummergengen} directly since for $L = \mathcal O_A$
it just states $0 = 0$.

However, it is still possible to make use of the proposition. We have to look
at all generalised Kummer varieties at the same time. Doing so we get the
following main result of this work:
\begin{thm}
  \label{thm:kummergen}
  Let $\phi: \Omega \to R$ be a complex genus with values in the $\Q$-algebra
  $R$. For every smooth projective surface $X$ with $\int_X c_1(X)^2 \neq
  0$,
  \begin{equation}
    \phi(K) = \frac 1{c_1(X)^2} \left(z \frac{\diff}{\diff z}\right)^2
    \ln \frac{\phi_1(H_X) \phi_{-1}(H_X)}{\phi(H_X)^2}.
  \end{equation}
\end{thm}

\begin{proof}
  Let $L$ be any line bundle on $A$. We have
  \begin{multline}
    c_1(L)^2 \sum_{n = 1}^\infty \phi(A^{[[n]]}, L^{[[n]]}) z^n
    = 2 \sum_{n = 1}^\infty n^2 \phi(A^{[n]}, L_n) z^n
    = 2 \left(z \frac{\diff}{\diff z}\right)^2 \phi(H_{A, L})
    \\
    = 2 \left(z \frac{\diff}{\diff z}\right)^2
    \exp\left(c_1(L)^2 A_\phi + c_1(L) c_1(A) B_\phi + c_1(A)^2 C_\phi
      + c_2(A) D_\phi \right)
    \\
    = 2 \left(z \frac{\diff}{\diff z}\right)^2
    \exp\left(c_1(L)^2 A_\phi\right)
    = 2 \left(z \frac{\diff}{\diff z}\right)^2 c_1(L)^2 A_\phi + \mathrm
    O\left((c_1(L)^2)^2\right),
  \end{multline}
  which together with proposition~\ref{prop:aphi} proves the theorem, since
  there are line bundles on $A$ with $c_1(L) \neq 0$.
\end{proof}

\begin{rmk}
  Of course, everything still holds true if we replace the abelian surface $A$
  from which we constructed the generalised Kummer varieties, by an arbitrary
  complex torus of dimension two.
\end{rmk}

\newpage

\appendix

\section{The Chern numbers of the generalised Kummer varieties of dimension up
  to 14}

We have used \ref{thm:kummergen} to compute all Chern numbers of the
generalised Kummer varieties of dimension up to fourteen. Our results are
as follows:
\begin{center}
  \begin{tabular}[t]{|c|r|}
    \hline
    Chern number & Evaluated on $A^{[[*]]}$ \\
    \hline
    $c_{2}$ & $24$ \\
    \hline
    $c_{2}^2$ & $756$ \\
    $c_{4}$ & $108$ \\
    \hline
    $c_{2}^3$ & $30208$ \\
    $c_{2}c_{4}$ & $6784$ \\
    $c_{6}$ & $448$ \\
    \hline
    $c_{2}^4$ & $1470000$ \\
    $c_{2}^2c_{4}$ & $405000$ \\
    $c_{4}^2$ & $111750$ \\
    $c_{2}c_{6}$ & $37500$ \\
    $c_{8}$ & $750$ \\
    \hline
    $c_{2}^5$ & $84478464$ \\
    $c_{2}^3c_{4}$ & $26220672$ \\
    $c_{2}c_{4}^2$ & $8141472$ \\
    $c_{2}^2c_{6}$ & $3141504$ \\
    $c_{4}c_{6}$ & $979776$ \\
    $c_{2}c_{8}$ & $142560$ \\
    $c_{10}$ & $2592$ \\
    \hline
    $c_{2}^6$ & $5603050432$ \\
    $c_{2}^4c_{4}$ & $1881462016$ \\
    $c_{2}^2c_{4}^2$ & $631808744$ \\
    $c_{4}^3$ & $212190776$ \\
    $c_{2}^3c_{6}$ & $268796752$ \\
    $c_{2}c_{4}c_{6}$ & $90412056$ \\
    $c_{6}^2$ & $12976376$ \\
    $c_{2}^2c_{8}$ & $17075912$ \\
    $c_{4}c_{8}$ & $5762400$ \\
    $c_{2}c_{10}$ & $441784$ \\
    $c_{12}$ & $2744$ \\
    \hline
  \end{tabular}
  \begin{tabular}[t]{|c|r|}
    \hline
    Chern number & Evaluated on $A^{[[*]]}$ \\
    \hline
    $c_{2}^7$ & $421414305792$ \\
    $c_{2}^5c_{4}$ & $149664301056$ \\
    $c_{2}^3c_{4}^2$ & $53149827072$ \\
    $c_{2}c_{4}^3$ & $18874417152$ \\
    $c_{2}^4c_{6}$ & $24230756352$ \\
    $c_{2}^2c_{4}c_{6}$ & $8610545664$ \\
    $c_{4}^2c_{6}$ & $3059945472$ \\
    $c_{2}c_{6}^2$ & $1397121024$ \\
    $c_{2}^3c_{8}$ & $1914077184$ \\
    $c_{2}c_{4}c_{8}$ & $681332736$ \\
    $c_{6}c_{8}$ & $110853120$ \\
    $c_{2}^2c_{10}$ & $71909376$ \\
    $c_{4}c_{10}$ & $25700352$ \\
    $c_{2}c_{12}$ & $1198080$ \\
    $c_{14}$ & $7680$ \\
    \hline
  \end{tabular}
\end{center}

It is remarkable fact that all Chern numbers of the varieties $A^{[[n]]}$ with
$n \leq 8$ are positive and divisible by $n^3$. As the known Chern numbers of
Hilbert schemes of points on K3 surfaces are also positive one can wonder if,
given a compact Hyperk\"ahler manifold $X$, all Chern numbers of $X$ are
positive.
 
\bibliographystyle{plain}
\bibliography{mybib}

\begin{thebibliography}{10}

\bibitem{beauville83}
Arnaud Beauville.
\newblock Vari\'et\'es {K}\"ahleriennes dont la premi\`ere classe de {C}hern
  est nulle.
\newblock {\em J. Differential Geom.}, 18(4):755--782, 1983.

\bibitem{britze01}
Michael Britze and Marc~A. Nieper.
\newblock Hirzebruch-{R}iemann-{R}och formulae on irreducible symplectic
  {K}\"ahler manifolds.
\newblock {\em arXiv:math.AG/0101062}.

\bibitem{dijkgraaf97}
Robert Dijkgraaf, Gregory Moore, Erik Verlinde, and Herman Verlinde.
\newblock Elliptic genera of symmetric products and second quantized strings.
\newblock {\em Comm. Math. Phys.}, 185:197--201, 1997.

\bibitem{lehn01}
Geir Ellingsrud, Lothar G\"ottsche, and Manfred Lehn.
\newblock On the cobordism class of the {Hilbert} scheme of a surface.
\newblock {\em J. Algebraic Geom.}, 10(1):81--100, 2001.

\bibitem{ellings87}
Geir Ellingsrud and Stein~Arild Str\o{}mme.
\newblock On the homology of the {H}ilbert scheme of points in the plane.
\newblock {\em Invent. Math.}, 87:343--352, 1987.

\bibitem{ellings96}
Geir Ellingsrud and Stein~Arild Str\o{}mme.
\newblock {B}ott's formula and enumerative geometry.
\newblock {\em J. Amer. Math. Soc}, 9:175--193, 1996.

\bibitem{fogarty68}
John Fogarty.
\newblock Algebraic families on an algebraic surface.
\newblock {\em Amer. J. Math}, 90:511--521, 1968.

\bibitem{goettsche93}
Lothar G\"ottsche and Wolfgang Soergel.
\newblock Perverse sheaves and the cohomology of {Hilbert} schemes of smooth
  algebraic surfaces.
\newblock {\em Math. Ann.}, 296(2):235--245, 1993.

\bibitem{hl97}
Daniel Huybrechts and Manfred Lehn.
\newblock {\em The {G}eometry of {M}oduli {S}paces of {S}heaves}, volume E31 of
  {\em Aspects of Mathematics}.
\newblock Friedr. Vieweg \and Sohn Verlagsgesellschaft mbH,
  Braunschweig/Wiesbaden, 1997.

\bibitem{lehn99}
Manfred Lehn.
\newblock {C}hern classes of tautological sheaves on {H}ilbert schemes of
  points on surfaces.
\newblock {\em Invent. Math.}, 136:157--207, 1999.

\bibitem{milnor60}
John Milnor.
\newblock On the cobordism ring ${\Omega}^*$ and a complex analogue.
\newblock {\em Amer. J. Math.}, 82:505--521, 1960.

\bibitem{mumford70}
David Mumford.
\newblock {\em Abelian varieties}.
\newblock Published for the Tata Institute of Fundamental Research, Bombay,
  1970.
\newblock Tata Institute of Fundamental Research Studies in Mathematics, No. 5.

\bibitem{sawon99}
Justin Sawon.
\newblock {\em {Rozansky-Witten} invariants of hyperk\"ahler manifolds}.
\newblock PhD thesis, University of Cambridge, October 1999.

\end{thebibliography}

\end{document}